# MULTIPLE INTELLIGENCES AND QUOTIENT SPACES

**ABSTRACT**

The Multiple Intelligence Theory (MI) is one of the models that study and describe the cognitive abilities of an individual. In [7] is presented a referential system which allows to identify the Multiple Intelligences of the students of a course and to classify the level of development of such Intelligences. Following this tendency, the purpose of this paper is to describe the model of Multiple Intelligences as a quotient space, and also to study the Multiple Intelligences of an individual in terms of this new mathematical representation.

## INTRODUCTION

If we accepted as definition of intelligence to the set of an individual's capacities, which allow him to solve daily problems, to generate new problems and to create products and/or to offer services inside the cultural environment in which such individual lives, then it is possible to associate to each person -according to the Cognitive Theory- at least eight intelligences, or eight cognitive capacities. These intelligences work together, although as semi-independent cognitive entities, which usually do not have a level of uniform development in the individual (i.e. ones can be developed more than others). Furthermore, the cultures and segments of the society make different emphases in them.

Among the models that study and describe the cognitive abilities of an individual we can find the so-called Multiple Intelligences Theory (MI), created by H. Gardner, (1983, [4]) -a neuropsychologist and educator of the Graduate School of Education of the Harvard University-. Originally, [4] was only aimed to psychologists, but it found fertile land between educators and mathematicians. The most important contribution of the Multiple Intelligences Theory to Education is that it allows to the educators to expand their repertoire of methods, tools and strategies beyond those that are frequently used in the classrooms.

Gardner begins with the existence of many and different intellectual faculties, or competencies. Each faculty or competence describes its own history of development. In last years, Neurobiology has indicated the presence of some areas or regions in the brain that correspond, to certain extent, to some types of cognition; and such facts imply a neural organization that might correspond to the notion of different ways of information processing.

It is necessary to highlight, that it does not exist, and will never be able to exist, a single, irrefutable and accepted list of human intelligences. There will never be a master list of 3; 7 or 100 intelligences which could be guaranteed by researchers. It might be possible that a decisive theory of human intelligences could never be completed. But the reason to continue researching about this subject is the necessity we have for better classifications of human intelligences than the ones we actually have nowadays; evidences of such necessity exist and they are product of scientific researches, transcultural observations and educational study.

In [7] is presented a referential system that allows to identify the Multiple Intelligences of the students from a course and to classify the level of development of such Intelligences. Taking into account this tendency, we will describe –via quotient spaces- such referential system. Therefore, the fundamental result of this paper allows identifying the Multiple Intelligences of a person with a concrete mathematical object (in this case a quotient space).

The outline of the paper continues as follows: Section 1 is a brief description of the Multiple Intelligences Theory, which is necessary for the development of this work. Section 2 presents a construction of a quotient space associated to the referential system given in [7].

## 1. BASIC FACTS ON MULTIPLE INTELLIGENCES

In this section we will show some definitions and characteristics related to the Multiple Intelligences Theory (we should not review the whole development of this Theory. The reader is referred to [1], [4] , [7] or [8] for a detailed presentation of these or more recent results).

We define the intelligence as the set of an individual's capacities, which allow him to solve daily problems, to generate new problems and to create products and/or to offer services inside the cultural environment in which such individual lives.

A cognitive model that studies and describes the intelligence regarding the previous definition is the H. Gardner's model. In order to describe his model, Gardner used the following facts:

1. The existence of people which suffer illnesses or accidents and specific areas of their brain are damaged.
2. The existence of people with mental retards, prodigies and other exceptional people type that show a very uneven profile of abilities and differences.

3. Each intelligence follows its own evolutionary patron whose development is different on one of each.
4. Each intelligence is product of the evolution (this evolution is appreciated so much in the human species as in other species).
5. There are discoveries of psychometric investigation that support the existence of diverse intelligences.
6. Psychological researches of experimental type affirm that each intelligence operates in separated way with respect to the other ones.
7. Identification of the development history of each intelligence -this fact occupies a primordial place in the education-.
8. Each intelligence possesses a symbols system own.

In this way, Gardner establishes in his MI model that each person possesses at least eight intelligences, these work together, although as semi-independent cognitive entities. We will give a brief description of them.

**Linguistic Intelligence:** capacity to use words of effective way, so much in oral form as a written. It includes ability in the use of the syntax, phonology, semantics and pragmatic functions. For example, writers, educators, lawyers and narrators possess this intelligence in a high development level.

**Logical-Mathematical Intelligence:** capacity to use numbers in an effective way, and to transform with dexterity, different reasoning chains. A person with a good development of logical-mathematical intelligence, highlights in the resolution of problems, to carry out complicated mathematics calculations and logical reasoning, as scientists, engineers, economists, administrators, accountants and others.

**Spatial Intelligence:** capacity to perceive the visual and spatial world, and to transform or to recognize its elements. This intelligence includes several informal abilities, such as: ability to use the imagination and then transform it, ability to visualize colors, lines, shapes and figures and others, ability to produce graphic likeness of spatial information and the orientation ability. Some professionals that possess this intelligence in a high level of development are the painters, designers, architects, astronauts, mathematicians and engineers.

**Bodily-Kinaesthetic Intelligence:** capacity to use the body in different forms and to work cleverly with objects. This intelligence possesses specific physical abilities as the coordination, the balance, the dexterity, the force, the flexibility and the speed. It is manifested in sportsmen, dancers, scenic and plastic artists and others.

**Musical Intelligence:** capacity to perceive, to discriminate, to transform and to express musical shapes. It includes medullary abilities as the tone, the rhythm and the ringer. The professionals that evidence this intelligence are: composers, interpreters, musicians, musical educators and others.

**Interpersonal Intelligence:** capacity to perceive, to understand and to distinguish the different states of encouragement, intentions, reasons and other people's feelings through the communication. High level of development in this intelligence is observed in psychologists, sociologists, educators, journalists, doctors and others.

**Intrapersonal Intelligence:** capacity to build a precise perception with respect to itself, and to organize and to direct your own life. This intelligence includes abilities as knowing your own ideas, dexterities and the personal goals. It is manifested in theologians, philosophers and in capable persons to recognize their encouragement state and feelings.

**Naturalist Intelligence:** capacity to distinguish, to classify and to use the numerous species of the flora and fauna in natural environments. It includes abilities to understand behaviors, necessities and characteristics of the animals and plants; also, to experience, to meditate and to question about our environment. Great part of these components is in the biologists, zoologists, agronomists, geographers, astronomers, etc.

It is well-known that each description model of the intelligence requires a test, which provides mainly the development level of intelligence for one or several persons. This test consists on the application of a questions series, tasks, stimuli, situations and others, which help to determine the capacity of person knowledge and to compare it with other persons.

There are many intelligence test, and many classifications for them. Such classifications depend of certain parameters as method, purpose, modality of the application, behavior area, type of problem and demands done to the person.

T. Armstrong (2000, [1]) proposed a MI-test for discovering and encouraging the children's Multiple Intelligences. For this, he separates certain cognitive abilities of children -obtained in terms of the observation- in blocks of MI. The ideas in [1] can be used to give a MI-test which identifies the level of Multiple Intelligences of the teens (see [7]).

The MI-test in [7] considers certain cognitive abilities which are presented in some MI, and the response given by the person are associated with those cognitive abilities of MI-test that the person recognizes to have, or with those that the person is identified.

Let a test of MI we call us Spider Web System (SWS) of a person to graphic representation obtained from such test computing the maximal number of abilities that the person recognizes to have in each Intelligence.

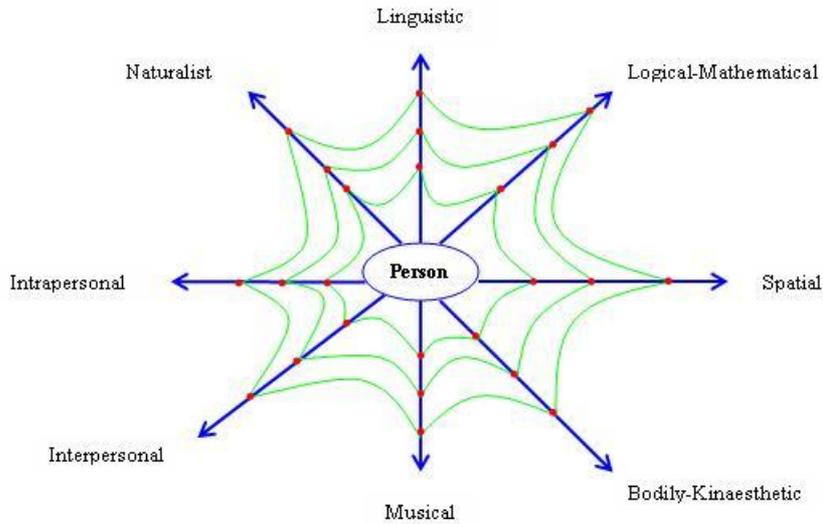

Figure 1. A person's ideal SWS.

The points represent the person's cognitive ability in each intelligence and the threads the interaction between two abilities corresponding to different intelligences.

Of this representation, we can deduce that there is not a patron of conditions that the person should gather to be intelligent in a particular area, so that the best way of diagnosing the intelligence is the observation.

Now, suppose that a person after having filled the MI-test obtained the following results, according to the sentences that he selected,

Linguistic Intelligence = 5          Logical-Mathematical Intelligence = 2
Spatial Intelligence = 6             Kinaesthetic-Bodily Intelligence = 7
Musical Intelligence = 4             Interpersonal Intelligence = 3
Intrapersonal Intelligence = 2       Naturalist Intelligence = 8

then the SWS associated will have the following form

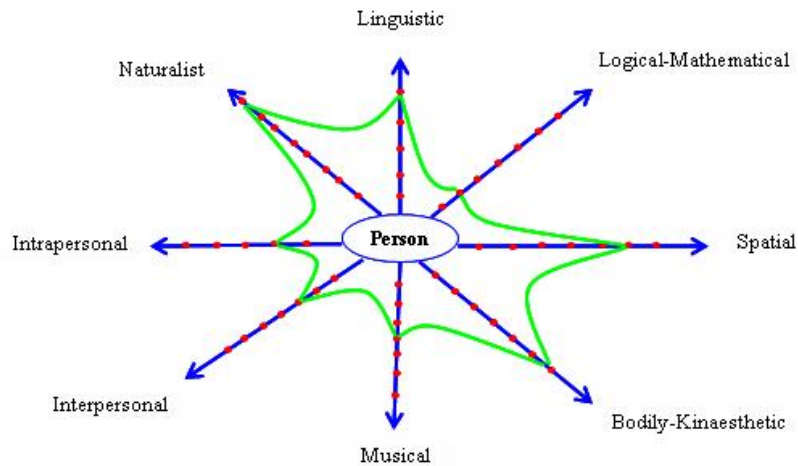

Figure 2. SWS corresponding to the test developed by an individual.

This example shows that we can try to classify two or more individuals in terms of their SWS respective, although these are not felt identified with the same sentences. Furthermore, we can group them so that the lack of person abilities is compensated by the abilities of others. It is to say, we can form groups where a person is less intelligent musically, but another balances this lack, being musically very intelligent. And this it is basically the objective of the test: to classify the persons in Homogeneous Intelligent Groups (HIG), in which the abilities lacking of a member are balanced with those of another of the same group. Then we would obtain work groups with SWS as the following

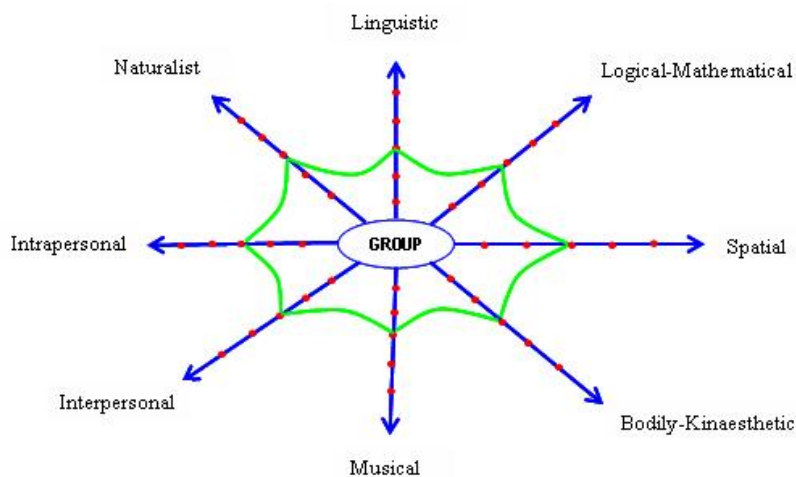

Figure 3. SWS of an HIG.

Notice that the MI-test can be adapted for persons of any age; since such adaptation depends of the present abilities in the test.

## 2. THE QUOTIENT SPACE OF THE MULTIPLE INTELLIGENCES

In this section we will present the main result of this paper, the quotient space associated to the reference system given in [7].

Let us consider $X$ the set of all cognitive abilities of a person -latent or no, developed or no-, it is clear that each MI is a subset of $X$. We denote by

$C_1 =$ Linguistic Intelligence  $\qquad C_2 =$ Logical-Mathematical Intelligence

$C_3 =$ Spatial Intelligence  $\qquad C_4 =$ Bodily- Kinaesthetic Intelligence

$C_5 =$ Musical Intelligence  $\qquad C_6 =$ Interpersonal Intelligence

$C_7 =$ Intrapersonal Intelligence  $\qquad C_8 =$ Naturalist Intelligence

Since a cognitive abilities $x$ can belong to several Multiple Intelligences, has sense to consider the following subset of $X$, given by $V := \{x \in X : x$ belongs to more of a MI$\}$.

Now, we will use to $V$ to construct a partition associated with the set of all MI, as follows.

For each $x \in V$ there are exist $i_1,...,i_k \in \{1,...,8\}$ such that $x \in C_{i_1} \cap ... \cap C_{i_k}$.

We can take $i = \min\{i_j : 1 \leq j \leq k\}$, it clear that $x \in \tilde{C_i} := C_i$.

Now, let us consider the following reduced intelligences $\tilde{C}_{i_j} = C_{i_j} | \{x\}$ for all $i_j \neq i$.

Then, by construction, the set $Y = \bigcup_{j=1}^{8} \tilde{C}_j$ is a partition for the set of all MI and it is well-know that every partition of a set induces a equivalence relation (see for example, [3]).

We say that $y \sim y'$ on $Y$, if and only if, there is exists $j \in \{1,...,8\}$ such that $y, y' \in \tilde{C}_j$.

In terms of this reference system, we can illustrate the previous equivalences classes as follows

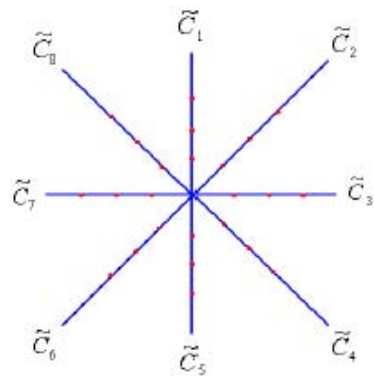

Figure 4. Graphic representation of the sets $\tilde{C}_j$.

In previous figure, the points on each axis represent the cognitive abilities corresponding to each set $\tilde{C}_j; j \in \{1,...,8\}$ and the axis represents the relationship between the cognitive abilities in the set.

On the other hand, if the Multiple Intelligences Theory emphasizes that, cognitive abilities corresponding to different intelligence, operate in group to solve certain problems -just as it is illustrated in the Figure 1-, then we can relate this operation when we take $Y$ as the union of the sets $\tilde{C}_j$. In such sense, graphic representation of $Y$ is:

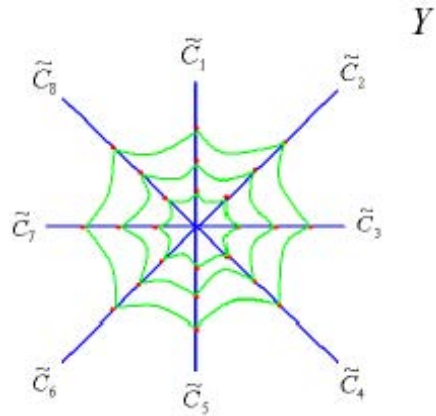

Figure 5. Graphic representation of the set $Y$.

Therefore, we have identified to the reference system given in [7] with a space quotient.